\documentclass[12pt]{article}
\usepackage{amsmath, amssymb, amsthm, bbm}
\usepackage{graphicx, psfrag, epsf, subcaption}
\usepackage{enumerate, algorithm}
\usepackage{natbib, hyperref}
\usepackage{booktabs, threeparttable}
\usepackage[title]{appendix}
\usepackage{url} % not crucial - just used below for the URL

\newtheorem{prop}{Proposition}
\newtheorem{them}{Theorem}

\theoremstyle{definition}
\theoremstyle{definition}
\theoremstyle{definition}\newtheorem{assu}{Assumption}
\theoremstyle{definition}\newtheorem*{remk}{Remark}
\graphicspath{{figs/}}
\hypersetup{colorlinks, linkcolor=blue, citecolor=blue, anchorcolor=blue, urlcolor=blue}

\begin{document}

\def\spacingset#1{\renewcommand{\baselinestretch}%
{#1}\small\normalsize} \spacingset{1}

%%%%%%%%%%%%%%%%%%%%%%%%%%%%%%%%%%%%%%%%%%%%%%%%%%%%%%%%%%%%%%%%%%%%%%%%%%%%%%

\title{\bf A factor-adjusted multiple testing of \\ general alternatives}
\author{Mengkun Du\thanks{Email: \url{dumengkun@pku.edu.cn}}\hspace{.2cm}\\
    School of Mathematical Sciences, Peking University\\
    and \\
    Lan Wu\thanks{Lan Wu is the corresponding author (Email: \url{lwu@pku.edu.cn}) at School of Mathematical Sciences, Peking University, Beijing, 100871, China.}\hspace{.2cm}\\
    School of Mathematical Sciences, Peking University}
\maketitle

\begin{abstract}
Factor-adjusted multiple testing is used for handling strong correlated tests. Since most of previous works control the false discovery rate under sparse alternatives, we develop a two-step method, namely the AdaFAT, for any true false proportion. In this paper, the proposed procedure is adjusted by latent factor loadings. Under the existence of explanatory variables, a uniform convergence rate of the estimated factor loadings is given. We also show that the power of AdaFAT goes to one along with the controlled false discovery rate. The performance of the proposed procedure is examined through simulations calibrated by China A-share market.
\end{abstract}

\noindent%
\vfill

\newpage
\spacingset{1} % DON'T change the spacing!

\section{Introduction}
\label{sec:int}

\cite{benjamini1995controlling} introduced the false discovery rate (FDR) as a criterion to qualify multiple testing procedures and showed that the proposed Benjamini and Hochberg (BH) procedure controls the FDR if p-values from true nulls are independent. Although this independence assumption can be relaxed to arbitrary dependency by \cite{benjamini2001the}, there is a conservative adjustment limiting the power of test. \cite{storey2004strong} further improved the BH procedure by allowing weak correlated p-values, but the FDR control is still challenged when p-values are strong correlated with each other.

To characterize such strong correlation, there are two approaches. One is to build models upon p-values or the corresponding z-scores that was suggested by \cite{efron2007correlation} and \cite{sun2009large}. Another is to introduce surrogate variables at the level of original data that was suggested by \cite{leek2008a} and \cite{friguet2009a}. In particular, \cite{friguet2009a} considered the following linear model:
\begin{equation}\label{equ:model-i}
  y_{ij}=\alpha_j+\beta_j'x_i+\gamma_j'z_i+\varepsilon_{ij}\quad\text{for}\quad i\leq n,\ j\leq m,
\end{equation}
where $n$ is the number of observations, $m$ is the number of tests, and $z_i$ refers to the $i$th realization of latent factors associated with testing dependency. The FDR control of $\alpha_j$ and $\beta_j$ can be effectively improved after adjusting by the latent factors. As a result, factor-adjusted multiple testing procedures based on equation (\ref{equ:model-i}) are broadly discussed in recent literatures.

\cite{fan2012estimating} provided an estimator to approximate the false discovery proportion. \cite{wang2017confounder} introduced a confounder between $x_i$ and $z_i$. \cite{fan2019farmtest:} proposed a robust procedure with Huber regression. \cite{Lan2019A} consistently estimated the number of latent factors. All preceding works are under the assumption of sparse alternatives, i.e. the proportion of true nulls is close to one for large $m$. Sparse condition is reasonable for microarray data, but for financial applications, \cite{barras2019a} inferred the proportion of non-zero $\alpha_j$'s and found that there are around $37.3\%$ mispriced assets in US market, which means that the condition of sparse $\alpha$ does not hold and should be removed under this circumstance.

In this paper, a two-step factor-adjusted multiple testing, namely the AdaFAT, is proposed to handle strong cross-sectional correlations. We let the original t-tests be a preprocessed procedure to estimate a subset of the true nulls, and use this subset to modify the factor-adjusted test statistics for general alternatives. We further extend the suitability of our method by allowing explanatory variables. When there is $x_i$ in equation (\ref{equ:model-i}), the estimator of $\gamma_j$ is derived from auto-correlated samples. We can improve theorems in \cite{bai2002determining} to obtain the estimator $\hat{\gamma}_j$ with a uniform convergence rate. Based on this uniform consistency of $\hat{\gamma}_j$, we show that the power of our method goes to one. Meanwhile, the FDR is controlled without the assumption of sparse alternatives. To our knowledge, this asymptotic result is obtained for the first time among relevant literatures.

The rest of paper is set as follows. We describe our method in section \ref{sec:des}. The analytical results are presented in section \ref{sec:mai}. Section \ref{sec:sim} outlines the numerical studies, and section \ref{sec:con} concludes the paper.

\begin{remk}
Throughout the paper, denote $\varepsilon_i\ i.i.d.\ N(\mu,\Sigma)$ as a series of independent and identically distributed normal random vectors. Let $\mathbbm{1}_n$ be an $n$-vector of ones. For any matrix $X$ of $n$ rows, let $X'$ be the transpose of $X$ and $Q(X)=I_n-X(X'X)^{-1}X'$ be the projection matrix orthogonal to $X$. Also, denote matrix norms as $\|X\|=\sqrt{\lambda_{max}(X'X)}$, $\|X\|_F=\sqrt{tr(X'X)}$, and $\|X\|_{\infty}=\max_{i}\sum_j|x_{ij}|$.
\end{remk}

\section{Model Description}
\label{sec:des}

\subsection{Model Setup}
Suppose that there are $m$ hypothesis tests
\begin{equation}\label{equ:hypotheses}
  H_{0j}:\alpha_j=0\quad\text{vs.}\quad H_{1j}:\alpha_j\neq0\quad\text{for}\quad j\leq m.
\end{equation}
Without loss of generality, to simplify the discussion, we rewrite equation (\ref{equ:model-i}) in matrix notation as follows,
\begin{equation}\label{equ:model}
  Y=\mathbbm{1}_n\,\alpha'+XB+Z\,\Gamma+E,
\end{equation}
where each row of $Y=(y_{1},\ldots,y_{n})'\in\mathbb{R}^{n\times m}$ represents an observation $y_i=(y_{i1},\ldots,y_{im})'$. Let $p$ and $q$ be fixed numbers of explanatory variables and latent factors, respectively, then $X=(x_1,\ldots,x_n)'\in\mathbb{R}^{n\times p}$ and $Z=(z_1,\ldots,z_n)'\in\mathbb{R}^{n\times q}$ consist of $x_i=(x_{i1},\ldots,x_{ip})'$ and $z_i=(z_{i1},\ldots,z_{iq})'$. The combination $Z\,\Gamma+E$ usually refers to model residuals. We assume the existence of unobservable $Z$ to quantify the cross-sectional correlation unexplained by $X$, therefore $E=(\varepsilon_{1},\ldots,\varepsilon_{n})'\in\mathbb{R}^{n\times m}$ can be seen as idiosyncratic component of weakly correlated vector $\varepsilon_{i}=(\varepsilon_{i1},\ldots,\varepsilon_{im})'$. Model parameters are of either constants or random variables independent of $\{x_i,z_i,\varepsilon_i\}_{i\leq n}$. We denote $\alpha=(\alpha_1,\ldots,\alpha_m)'\in\mathbb{R}^{m}$ as the testing interval, and $B=(\beta_1,\ldots,\beta_m)\in\mathbb{R}^{p\times m}$ and $\Gamma=(\gamma_1,\ldots,\gamma_m)\in\mathbb{R}^{q\times m}$ as factor loadings of explanatory variables and latent factors, respectively.

We first illustrate the influence of latent factors. If $Z$ is treated as part of residuals, then the original t-tests along cross-sections can be formed as
\begin{equation}\label{test:ori}
  T_{ori}=\widehat{\Lambda}_{z\varepsilon}^{-1/2}\,\vartheta\in\mathbb{R}^m,
\end{equation}
where each entry $t_{ori,j}$ is equal to the t-score towards $H_{0j}$, diagonal matrix $\widehat{\Lambda}_{z\varepsilon}$ is from the least square estimator of $var(\Gamma'z_i+\varepsilon_i)$, and
\begin{equation}\label{equ:reg}
  \vartheta:=\frac{Y'Q(X)\mathbbm{1}_n}{\sqrt{\mathbbm{1}_n'Q(X)\mathbbm{1}_n}}=:c_n\alpha+\Gamma'\zeta+\eta,
\end{equation}
where $c_n=\sqrt{\mathbbm{1}_n'Q(X)\mathbbm{1}_n}$, $\zeta=Z'Q(X)\mathbbm{1}_n/c_n$, and $\eta=E'Q(X)\mathbbm{1}_n/c_n$. Under the existence of $Z$, off-diagonal entries of $cov(T_{ori})$ are uniformly far from zero, which leads to strong correlations among the corresponding p-values. In this case, we temporarily assume that unknown parameters are observable and consider the oracle procedure
\begin{equation}\label{test:ora}
  T_{ora}=\Lambda_{\varepsilon}^{-1/2}\,\left(\vartheta-\Gamma'\zeta\right),
\end{equation}
which is standardized by the diagonal matrix of idiosyncratic variances $\Lambda_{\varepsilon}$. If idiosyncratic errors are weakly dependent, traditional multiple testing like the BH procedure based on $T_{ora}$ has the controlled FDR. Thus, we let $T_{ora}$ be a target and focus on the following factor-adjusted test statistics
\begin{equation}\label{test:adj}
  T_{adj}=\widehat{\Lambda}_{\varepsilon}^{-1/2}\,\left(\vartheta-\widehat{\Gamma}'\hat{\zeta}\right).
\end{equation}

\subsection{Methodology}

Estimators in equation (\ref{test:adj}) are given in this subsection for getting a consistent $T_{adj}$. To begin with, denote $\widetilde{X}=(\mathbbm{1}_n,X)$, and we use the principal component decomposition of
\begin{equation*}
  \widetilde{Y}:=Q(\widetilde{X})Y=Q(\widetilde{X})Z+Q(\widetilde{X})E=:\widetilde{Z}+\widetilde{E},
\end{equation*}
to obtain the estimators $\widehat{\Gamma}$ and $\widehat{\Lambda}_{\varepsilon}$. Since $\widehat{Z}$ as the eigenvectors of $\widetilde{Y}\widetilde{Y}'$ is orthogonal to $\widetilde{X}$, we cannot directly recover $\zeta$ from $\widehat{Z}$. If we substitute $\widehat{\Gamma}$ into equation (\ref{equ:reg}) and infer $\zeta$ by regressions, then there is a bias of non-zero $\alpha$. Therefore, a key issue for factor-adjusted multiple testing of general alternatives is to find a proper estimator $\hat{\zeta}$. We first apply a cross-sectional regression to equation (\ref{equ:reg}) and propose the steps below as a base algorithm.

\begin{algorithm}
\label{alg-fat}
\caption{FAT-DW}
\begin{enumerate}[(S1)]
  \item Consider the $n\times n$ matrix $\,\widetilde{Y}\widetilde{Y}'$, take its eigenvectors $(\xi_1, \ldots, \xi_k)$ corresponding to the largest $k$ eigenvalues, and define
  \begin{equation}\label{equ:ZGamma}
  \widehat{Z}_k=\sqrt{n}(\xi_1, \ldots, \xi_k)\quad\text{and}\quad\widehat{\Gamma}_k=\frac{1}{n}\widehat{Z}_k'Y;
  \end{equation}
  \item Let $\kappa$ be an upper bound such that $\kappa>q$ and $g(m,n)$ be any penalty function satisfying
  \begin{equation*}
    g(m,n)\to0\quad\text{and}\quad \min\{m,n\}g(m,n)\to\infty,
  \end{equation*}
  as $min\{m,n\}\to\infty$, minimize the following information criterion
  \begin{equation}\label{equ:ic}
  IC\left(k\right)=
  \log\left[\frac{1}{mn}\sum_{j=1}^{m}\sum_{i=1}^{n}(\widetilde{y}_{ij}-\hat{\gamma}_j'\hat{z}_i)^2\right]+kg(m,n),
  \end{equation}
  to obtain $\hat{q}$, and denote $\widehat{Z}=\widehat{Z}_{\hat{q}}$ and $\widehat{\Gamma}=\widehat{\Gamma}_{\hat{q}}\,;$
  \item Let $\widehat{E}=\widetilde{Y}-\widehat{Z}\widehat{\Gamma}$ be the estimated residuals, and obtain
  \begin{equation}\label{equ:zeta}
  \hat{\zeta}=\left[\widehat{\Gamma}Q(\mathbbm{1}_m)\widehat{\Gamma}'\right]^{-1}
  \widehat{\Gamma}Q(\mathbbm{1}_m)\vartheta\quad\text{and}\quad
  \widehat{\Lambda}_{\varepsilon}=diag\left(\frac{1}{n}\widehat{E}'\widehat{E}\right).
  \end{equation}
\end{enumerate}
\end{algorithm}

\begin{remk}
Multiple testing is usually applied to cases when $m\gg n$. The FAT-DW decomposes from an $n\times n$ matrix and is time-efficient in computation. Equation (\ref{equ:ic}) is from \cite{bai2002determining}. In equation (\ref{equ:zeta}), centralized operator $Q(\mathbbm{1}_m)$ is crucial when $\alpha$ has a significant mean. The diagonal $\hat{\Lambda}_{\varepsilon}$ indicates that we construct every test statistics by its own standard deviation, which is not equivalent to the independent idiosyncratic errors. Condition (C\ref{con:5}) in the next section specifies the cross-sectional correlation of $\varepsilon_i$.
\end{remk}

From equation (\ref{equ:zeta}), we deduce that the FAT-DW has an error term
\begin{equation*}
  \max_{j\in\mathcal{I}_1}\left|\gamma_j'(\Gamma Q(\mathbbm{1}_m)\Gamma')^{-1}\Gamma Q(\mathbbm{1}_m)\alpha\right|,
\end{equation*}
and is close to the oracle procedure when $\alpha$ is uncorrelated with each row of $\Gamma$. For more general cases of alternatives, however, this error decays only if $\alpha$ approaches to zero. To solve this problem, we find that although the original t-tests cannot guarantee a controlled FDR, it has the power of test that converges to one as sample size diverges. Thus, we use the original t-tests as a preprocessed procedure to estimate a subset of the true nulls, and $\zeta$ can be recovered from this subset when $m$ is sufficiently large. In practice, we propose an iterative algorithm of the FAT-DW to augment the test stability and reduce its bias to the oracle procedure.

\begin{algorithm}
\label{alg-adafat}
\caption{AdaFAT}
\begin{enumerate}[(S1)]
  \item Let $\,\hat{\mathcal{I}}_{ori,1}$ be the set of rejections from the original t-tests, and denote the initial subset of the true nulls as $\,\hat{\mathcal{I}}_0=\{1,\ldots,m\}\,\backslash\,\hat{\mathcal{I}}_{ori,1}$;
  \item Let $\widehat{\Gamma}_0$ and $\vartheta_0$ be matrices of the rows of $\,\widehat{\Gamma}$ and $\vartheta$ in $\hat{\mathcal{I}}_0$, respectively, let $\hat{m}_0$ be the cardinality of $\,\hat{\mathcal{I}}_0$, and modify the FAT-DW by renewing
  \begin{equation}\label{equ:zeta-k}
  \hat{\zeta}=\left[\widehat{\Gamma}_0Q(\mathbbm{1}_{\hat{m}_0})\widehat{\Gamma}_0'\right]^{-1}
  \widehat{\Gamma}_0Q(\mathbbm{1}_{\hat{m}_0})\vartheta_0;
  \end{equation}
  \item Run multiple testing with the modified $T_{adj}$ to get the rejection set $\,\hat{\mathcal{I}}_{adj,1}$, and update $\hat{\mathcal{I}}_0$ to $\hat{\mathcal{I}}_0\,\backslash\,\hat{\mathcal{I}}_{adj,1}$;
  \item Repeat (S2) and (S3) until $\,\hat{\mathcal{I}}_{adj,1}$ is stable.
\end{enumerate}
\end{algorithm}

\begin{remk}
The estimators $\,\widehat{\Gamma}$ and $\widehat{\Lambda}_{\varepsilon}$ are estimated from the complete set $\{1,\ldots,m\}$ and remain unchanged during iterations. The shrinking set $\,\hat{\mathcal{I}}_0$ is used for estimating $\zeta$, and the stable set of rejections $\,\hat{\mathcal{I}}_{adj,1}$ is the only output of interest.
\end{remk}

\subsection{Multiple Testing Framework}

The end of this section lists necessary notations of multiple testing procedures. Denote $t_{\cdot,j}$ as the $j$th entry of multiple test statistics and
\begin{equation}\label{equ:pval}
  p_{\cdot,j}=2\Phi(-|t_{\cdot,j}|),
\end{equation}
as the corresponding p-value towards double-side hypotheses. Let $\mathcal{I}_0$ be the set of true nulls having $m_0$ elements and $\mathcal{I}_1$ be the set of alternatives with cardinality $m_1$. For any threshold $t\in(0,1]$, denote
\begin{equation}
  \begin{aligned}
  V_{\cdot}(t)&=\#\{j\in\mathcal{I}_0:p_{\cdot,j}\leq t\},\\
  S_{\cdot}(t)&=\#\{j\in\mathcal{I}_1:p_{\cdot,j}\leq t\},\\
  R_{\cdot}(t)&=\#\{j\leq m:p_{\cdot,j}\leq t\},
  \end{aligned}
\end{equation}
as the number of false discoveries, the number of positive discoveries, and the number of rejections, respectively. Furthermore, we define the false discovery proportion and the power of test separately as follows,
\begin{equation}\label{equ:indi}
  \mathrm{FDP}_{\cdot}(t)=\frac{V_{\cdot}(t)}{\max\{R_{\cdot}(t),1\}}\quad\text{and}\quad
  \mathrm{POW}_{\cdot}(t)=\frac{S_{\cdot}(t)}{\max\{m_1,1\}}.
\end{equation}
Similar to type \uppercase\expandafter{\romannumeral1} and type \uppercase\expandafter{\romannumeral2} errors of single hypothesis tests, the false discovery proportion and the power are two key indicators that describe the superiority of multiple testing procedures. \cite{genovese2002operating} introduced the false non-discovery rate that is negatively correlated with the power of test in equation (\ref{equ:indi}). According to \cite{storey2004strong}, the false discovery rate can be inferred as follows,
\begin{equation}\label{equ:st}
\widehat{\mathrm{FDR}}_{\nu}(t)=\frac{\hat{\pi}_0(\nu)mt}{\max\{R_{\cdot}(t),1\}},
\end{equation}
where $\hat{\pi}_0(\nu)=[m-R_{\cdot}(\nu)]/[(1-\nu)m]$ is an estimated true null proportion with tuning parameter $\nu\in[0,1)$. If the estimator $\,\widehat{\mathrm{FDR}}_{\nu}(t)$ converges to $\mathrm{FDR}_{\cdot}(t)$ almost surely, then for any predetermined level $\tau>0$,
\begin{equation}\label{equ:tstar}
t^*=\sup\{t\in(0,1]:\widehat{\mathrm{FDR}}_{\nu}(t)\leq\tau\},
\end{equation}
will be an ideal threshold such that $\limsup\mathrm{FDR}_{\cdot}(t^*)\leq\tau$.

\section{Analytical Results}
\label{sec:mai}

\subsection{Model Assumptions}

Assume that there exist a $q\times q$ positive definite matrix $\Sigma_{\gamma}$ and constants $0<b_{min}\leq b_{max}<\infty$, $h<1$ satisfying
\begin{enumerate}[(C1)]
    \item $b_{min}\leq\lambda_{k}(\widetilde{X}'\widetilde{X}/n)\leq b_{max}$
          for $i\leq n,\,k\leq p+1$; \label{con:1}
    \item $(z_i',\varepsilon_i')\ i.i.d.\ N(0,diag\{I_q,\Sigma_{\varepsilon}\})$ that is also
          independent of $x_i$ for $i\leq n$;\label{con:2}
    \item $|\alpha_j|\geq b_{min}$ for $\mathcal{I}_1$, $\|Q(\mathbbm{1}_m)\alpha\|=O(m^h)$ and
          $\|\Gamma Q(\mathbbm{1}_m)\alpha\|=o(m/\sqrt{n})$;\label{con:3}
    \item $\|\gamma_j\|\leq b_{max}$ and
          $\lim_m\|\Gamma Q(\mathbbm{1}_m)\Gamma'/m-\Sigma_{\gamma}\|=0$ for $j\leq m$;\label{con:4}
    \item $\|\Sigma_{\varepsilon}\|_{\infty}\leq b_{max}$ and
          $\,b_{min}\leq\sigma_{\varepsilon,jj}\leq b_{max}$ for $j\leq m$;\label{con:5}
    \item $b_{min}\leq|\alpha_j|\leq b_{max}$ for $j\in\mathcal{I}_1$. \label{con:6}
  \end{enumerate}

The following assumptions ensure that the FAT-DW and the AdaFAT converges to the oracle procedure in large samples.

\begin{assu}\label{assu1}
Include condition (C\ref{con:1}) to condition (C\ref{con:5}).
\end{assu}

\begin{assu}\label{assu2}
Follow Assumption \ref{assu1} and replace (C\ref{con:3}) with (C\ref{con:6}).
\end{assu}

Conditions (C\ref{con:1}) and (C\ref{con:2}) are a trade-off between the existence of $x_i$ and the robustness of $(z_i',\varepsilon_i')$. If there is no explanatory variable, \cite{fan2018large} gives the possibility of replacing the normality assumption with elliptical distributions. If eigenvectors are estimated after screening out $x_i$, auto-correlated $\tilde{y}_i$ leads us to consider the normal innovations.

The lower bound and $\|Q(\mathbbm{1}_m)\alpha\|=O(m^h)$ in condition (C\ref{con:3}) are not stringent, which ensures the power and limits the dispersion of $\alpha$, respectively. $\|\Gamma Q(\mathbbm{1}_m)\alpha\|$ should be checked because $\|\Gamma Q(\mathbbm{1}_m)\alpha\|=o(m/\sqrt{n})$ is no longer fulfilled when $\alpha_j$ and $\gamma_j$ are cross-sectionally correlated. For significant $\|\Gamma Q(\mathbbm{1}_m)\alpha\|$, a modified condition is given in (C\ref{con:6}).

Conditions (C\ref{con:4}) and (C\ref{con:5}) are mainly set for the consistency of the estimated eigen-structure. The positive definite matrix assumes the pervasive factors of $z_i$ that has a global impact on $y_j$. See \cite{fan2019farmtest:} for more information of the pervasiveness condition. We assume $\|\Sigma_{\varepsilon}\|_{\infty}\leq b_{max}$ as in \cite{fan2013large} to constrain the sparsity of idiosyncratic components. This assumption can be relaxed to $\lambda_{max}(\Sigma_{\varepsilon})<b_{max}$ if the number of latent factors $q$ is known.

\subsection{Main Results}

We first show that the oracle procedure is consistent as follows. All proofs in this subsection are provided in appendices at the end of this paper.
\begin{prop}\label{pro:ora}
  Given $t^*$ in equation (\ref{equ:tstar}) with tuning parameter $\nu\in[0,1)$, suppose there is $\pi_0\in[0,1]$ such that $\lim_m m_0/m=\pi_0$ and $\log m=o(n)$, then as $min\{m_0, m_1, n\}\to\infty$, we have
  \begin{enumerate}[(a)]
    \item $\lim\sup\mathrm{FDR}_{ora}(t^*)\leq\tau$;
    \item $\lim\mathbb{E}\left[\mathrm{POW}_{ora}(t^*)\right]=1$;
  \end{enumerate}
  under Assumption \ref{assu1} or Assumption \ref{assu2}.
\end{prop}

Proposition \ref{pro:ora} holds for any true null proportion $\pi_0\in[0,1]$. If $\pi_0=1$, the BH procedure is good enough under this condition of sparse alternatives. If $\pi_0<1$, the multiple testing with $t^*$ is more powerful because $\hat{\pi}_0(\nu)$ broadens the rejection region of every single test. Next, we give asymptotics of the estimated eigen-structure in the FAT-DW and the AdaFAT.

\begin{prop}\label{pro:pca}
  The following consistency holds under Assumption \ref{assu1} or Assumption \ref{assu2}:
  \begin{enumerate}[(a)]
    \item $P(\hat{q}=q)\to1$ as $\min\{m,n\}\to\infty$;
    \item $\max_{j\leq m}\|\hat{\gamma}_j-H^{-1}\gamma_j\|=O_p(\sqrt{(\log m)/n})+O_p(1/\sqrt{m})$;
    \item $\max_{j\leq m}|\hat{\sigma}_{\varepsilon,jj}-\sigma_{\varepsilon,jj}|=
                            O_p(\sqrt{(\log m)/n})+O_p(1/\sqrt{m})$.
  \end{enumerate}
\end{prop}

In Proposition \ref{pro:pca}, $\widehat{\Gamma}$ is biased by a $q\times q$ matrix
\begin{equation}\label{equ:H}
  H=\frac{1}{mn}\Gamma\Gamma'\widetilde{Z}'\widehat{Z}D^{-1},
\end{equation}
where the diagonal matrix $D$ consists of eigenvalues such that $\widetilde{Y}\widetilde{Y}'\widehat{Z}/(mn)=\widehat{Z}D$. If we only apply the FAT-DW as well as the AdaFAT to factor-adjusted test statistics $T_{adj}$, then this identification issue does not affect our results because $T_{adj}$ is adjusted by
\begin{equation*}
  \widehat{\Gamma}'\hat{\zeta}=\widehat{\Gamma}'
  [\widehat{\Gamma}Q(1_m)\widehat{\Gamma}']^{-1}\widehat{\Gamma}Q(1_m)\vartheta,
\end{equation*}
that is unique for invertible $H$. Finally, the theorem below indicates that our method performs as close as the oracle procedure in large samples.

\begin{them}\label{the:fdp}
  Given $\,\log m=o(\sqrt{n})$ for any $t\in(0,1]$, as $min\{m, n\}\to\infty$, we have
  \begin{enumerate}[(a)]
    \item $|\mathrm{FDP}_{ora}(t)-\mathrm{FDP}_{adj}(t)|\to_{L_1}0$;
    \item $|\mathrm{POW}_{ora}(t)-\mathrm{POW}_{adj}(t)|\to_{L_1}0$,
  \end{enumerate}
  which holds for the FAT-DW under Assumption \ref{assu1} and for the AdaFAT under Assumption \ref{assu2}.
\end{them}

\section{Numerical Results}
\label{sec:sim}

The performance of the AdaFAT is examined in this section. We calibrate model parameters by weekly stock returns, and data used is from \emph{Wind Information Ltd.}, including daily prices and trading status of all stocks in China A-share market from January $2014$ to December $2018$.

There are $1634$ stocks being selected under a data-cleaning process of \cite{liu2019size}. Let market factor be the explanatory variable and $\alpha$ and $B$ be their ordinary least square estimators. We follow the suggestion of \cite{bai2003inferential} to denote the penalty function as $g(m,n)=\frac{m+n}{mn}\log\left(\frac{mn}{m+n}\right)$, and three latent factors are selected. We estimate $\Sigma_{\varepsilon}$ as the POET-estimator in \cite{fan2013large} and sample $\{x_i,z_i,\varepsilon_i\}$ from student's $t_{3}$ distribution for robustness concern.

\begin{figure}
    \centering
    \begin{subfigure}[b]{\textwidth}
        \includegraphics[width=\textwidth]{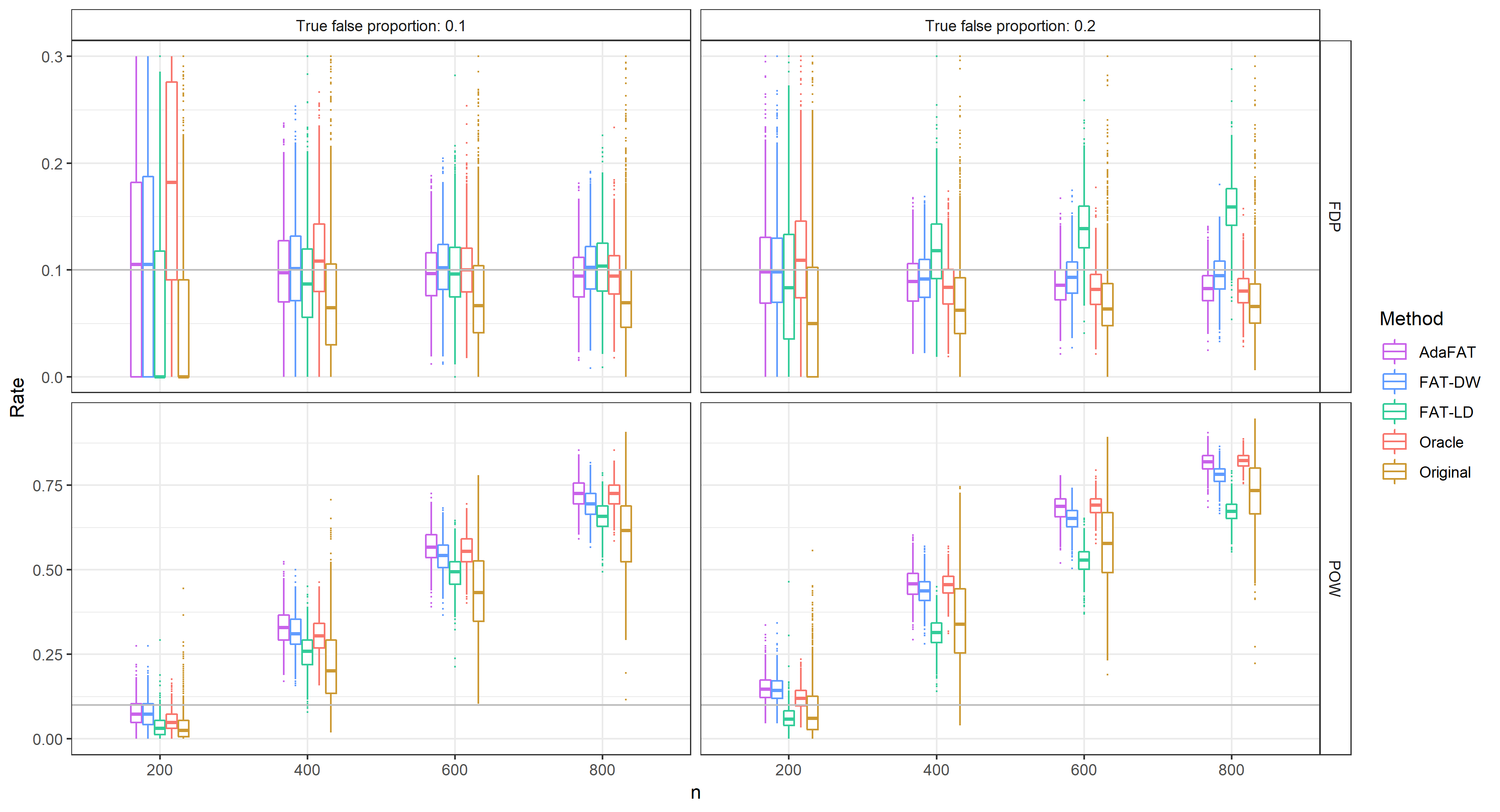}
        \caption{$\mu_z=0$.}
        \label{fig:simu:0}
    \end{subfigure}
    \vspace{0.1cm}

    \begin{subfigure}[b]{\textwidth}
        \includegraphics[width=\textwidth]{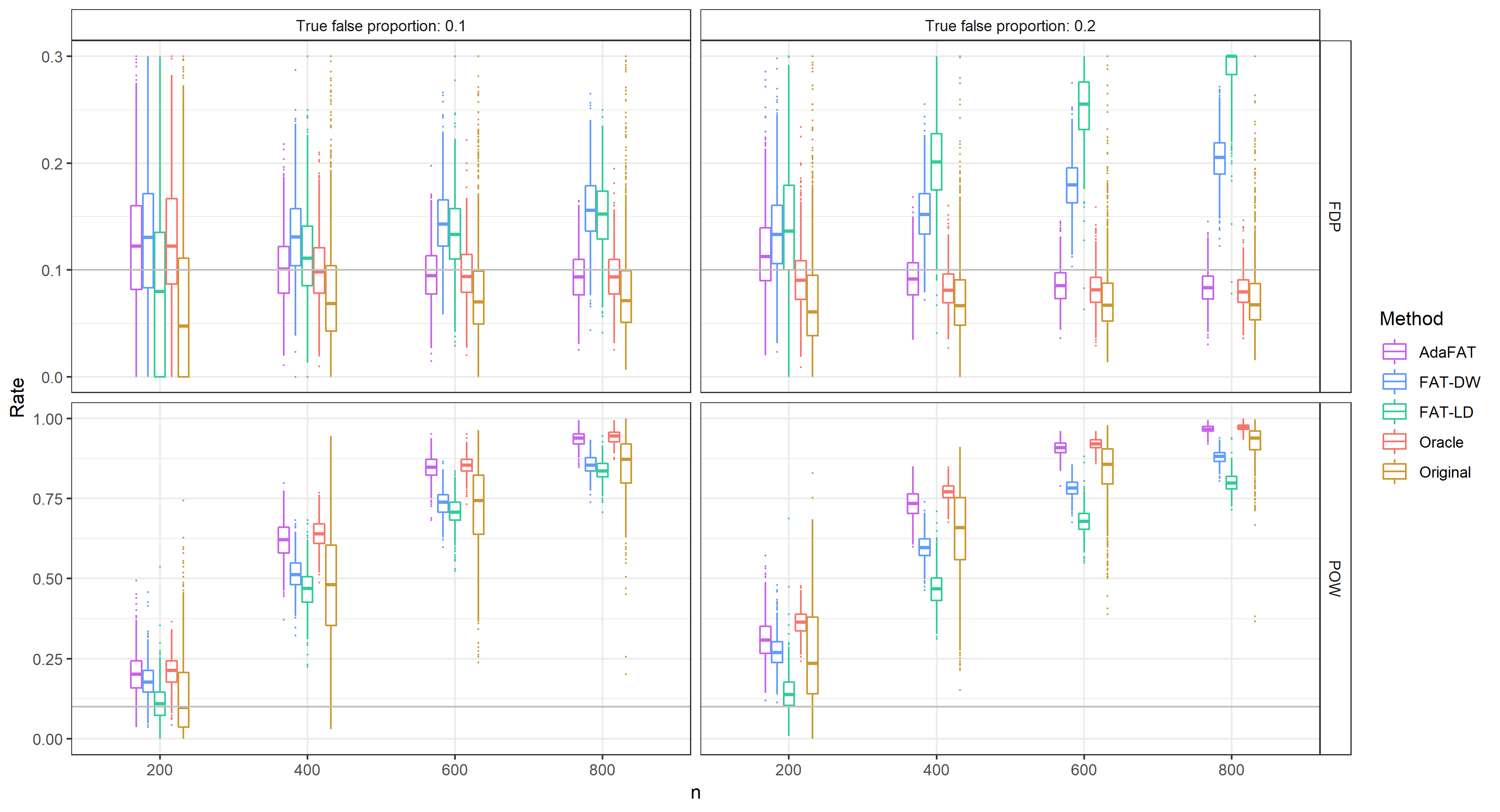}
        \caption{$\mu_z=\frac{1}{2}\mu_x\mathbbm{1}_3$.}
        \label{fig:simu:5}
    \end{subfigure}
    \caption{Boxplots of the simulated FDP and POW with $1000$ replications, where grey lines represent the FDR level, x-axis refers to the sample size and each column corresponds to a specific true false proportion $\pi_1$.}\label{fig:simu}
\end{figure}

Figure \ref{fig:simu} draws boxplots of $1000$ replicated results, where grey lines refer to the FDR level $\tau=0.1$. We compare our method with the original t-tests, the oracle procedure and \cite{Lan2019A}'s FAT-LD procedure. The difference in FAT-LD is that it decomposes eigenvectors from $Q(X)Y$, and $\|\alpha\|$ as an error term will bias $T_{adj}$ if the true false proportion $\pi_1$ is significant. We first set $\pi_1=0.1$ as in \cite{Lan2019A} to discuss the testing performance of sparse cases. Denote $|\alpha_j|/\sqrt{\sigma_{\varepsilon,jj}}$ as information ratio of the $j$th asset, and $\mathcal{I}_1$ consists of those of the largest $10\%$ information ratio. For non-sparse cases, we shift $\alpha$ to enlarge $\pi_1=0.2$ as shown in the second column of Figure \ref{fig:simu}. Also, we add an term for simulated samples as $\widetilde{\alpha} = \alpha + \Gamma'\mu_z$, in which the extra mean of latent factors $\mu_z$ enlarges the correlation between $\alpha$ and $\Gamma$. Figure \ref{fig:simu:0} refers to the original scenario of $\mu_z=0$, and we set each entry of $\mu_z$ to be $50\%$ of the mean of market factor $\mu_x$ in Figure \ref{fig:simu:5}.

The first column of Figure \ref{fig:simu} indicates that the oracle procedure converges in the way of Proportion \ref{pro:ora} with increasing sample size. On the contrary, the original t-tests has the biggest variance of FDP and POW, which implies the necessity of factor adjustment. In Figure \ref{fig:simu:0}, all factor-adjusted procedures are close to the oracle procedure when $\pi_1=0.1$, but when $\pi_1=0.2$, the second column shows that the FDR of FAT-LD is out of control. The FAT-DW performs more robust than the FAT-LD, and the AdaFAT is the closest one to the oracle procedure, which supports the statement of Theorem \ref{the:fdp}. In Figure \ref{fig:simu:5}, since there is a considerable $\mu_z$, non-sparse alternatives also bias the FAT-DW as sample size varies. Under this circumstance, the modified AdaFAT indicates its effectiveness with a lower FDR as well as a higher testing power than other factor-adjusted multiple testing procedures.

\section{Conclusion}
\label{sec:con}

In this paper, we propose the AdaFAT for strong correlated tests with general alternatives. As the cross-sectional dimension and sample size diverges, we show that the AdaFAT can find almost all non-nulls in probability with the controlled FDR. Simulations also support the suitability of AdaFAT. Regarding improvements for further research, the estimator of $\zeta$ should be considered with a more robust method to improve the testing efficiency especially for finite samples.

\section*{Acknowledgments}
\label{sec:ack}

Lan Wu's research was supported by Key Laboratory of Mathematical Economics and Quantitative Finance at Peking University.

\newpage
\begin{appendices}
\section{Proof of Proposition \ref{pro:ora}}
\label{app:ora}
\setcounter{equation}{0}
\renewcommand{\theequation}{A.\arabic{equation}}

For any constant $\delta>0$, the key point is to obtain the following equation
\begin{equation*}
    \frac{1}{m_1^2}\sum_{j,k\in\mathcal{I}_1}\left|cov\left(1_{\{p_{ora,j}\leq t\}},1_{\{p_{ora,k}\leq t\}}\right)\right|=O\left(\frac{1}{m_1^{\delta}}\right),
\end{equation*}
which is a sufficient condition of \cite{lyons1988strong}'s strong law of large numbers. \cite{fan2012estimating} deduce this condition under true nulls. Under alternatives, we have
\begin{equation*}
    \begin{aligned}
    &\left|cov\left(1_{\{p_{ora,j}\leq t\}},\,1_{\{p_{ora,k}\leq t\}}\right)\right|\\
    =&\left|P\left(|t_{ora,j}|\geq \psi_{t/2},|t_{ora,k}|\geq \psi_{t/2}\right)
      -P\left(|t_{ora,j}|\geq \psi_{t/2}\right)P\left(|t_{ora,k}|\geq \psi_{t/2}\right)\right|\\
    =&\left|P\left(|t_{ora,j}|< \psi_{t/2},|t_{ora,k}|< \psi_{t/2}\right)
      -P\left(|t_{ora,j}|< \psi_{t/2}\right)P\left(|t_{ora,k}|< \psi_{t/2}\right)\right|\\
    \leq&2\min\left\{P\left(|t_{ora,j}|< \psi_{t/2}\right),\,P\left(|t_{ora,k}|< \psi_{t/2}\right)\right\},
    \end{aligned}
\end{equation*}
where $\psi_x=\Phi^{-1}(1-x)$ is the upper quantile of standard normal distribution. Denote $\mu_j=c_n\alpha_j/\sqrt{\sigma_{\varepsilon,jj}}$ and $\eta_j^*=\eta_j/\sqrt{\sigma_{\varepsilon,jj}}$ such that
\begin{equation*}
  t_{ora,j}=\mu_j+\eta_j^*.
\end{equation*}
If $|t_{ora,j}|<\psi_{t/2}$, then $|\eta_j^*|>|\mu_j|-\psi_{t/2}$, which indicates that
\begin{equation*}
   \max_{j,k\in\mathcal{I}_1}\left|cov\left(1_{\{p_{ora,j}\leq t\}},\,1_{\{p_{ora,k}\leq t\}}\right)\right|
   \leq2P\left(|\eta_1^*|>\min_{j\in\mathcal{I}_1}|\mu_j|-\psi_{t/2}\right),
\end{equation*}
where $\eta_1^*\sim N(0,1)$ by definition. Also, we have
\begin{equation*}
  c_n\asymp\sqrt{n}\quad\text{and}\quad\min_{j\in\mathcal{I}_1}|\mu_j|\asymp\sqrt{n},
\end{equation*}
by conditions (C\ref{con:1}), (C\ref{con:3}) and (C\ref{con:6}). As $n\to\infty$, there is a $c_t>0$ such that
\begin{equation*}
   m_1^{\delta}\max_{j,k\in\mathcal{I}_1}\left|cov\left(1_{\{p_{ora,j}\leq t\}},\,1_{\{p_{ora,k}\leq t\}}\right)\right|
   \leq4m_1^{\delta}e^{-c_tn}\to0.
\end{equation*}
Thus, as $min\{m_0,m_1,n\}\to\infty$, we have strong law of large numbers below
\begin{equation*}
    \frac{V_{ora}(t)}{m_0}\to_{a.s.}t\quad\text{and}\quad \frac{S_{ora}(t)}{m_1}\to_{a.s.}1,
\end{equation*}
and the Theorem 4 in \cite{storey2004strong} completes the proof.

\section{Proof of Proposition \ref{pro:pca}}
\label{app:pca}
\setcounter{equation}{0}
\renewcommand{\theequation}{B.\arabic{equation}}

For any $k\geq 0$, we have the following identity
\begin{equation*}
  \widehat{Z}_k-\widetilde{Z}H_k=\frac{1}{mn}Q(\widetilde{X})\left(Z\Gamma E'+E\Gamma'Z'+EE'\right)\widehat{Z}_k,
\end{equation*}
where $H_k$ is defined as in \cite{bai2002determining}. Denote
\begin{equation*}
  M_X=\widetilde{X}(\widetilde{X}'\widetilde{X})^{-1}\widetilde{X}'
  \quad\text{and}\quad M_{ZE}=Z\Gamma E'+E\Gamma'Z'+EE',
\end{equation*}
for ease of discussion. The mean square error of $\widehat{Z}_k$ can be expressed as
\begin{equation}\label{equ:error:z}
  \begin{aligned}
  \frac{1}{n}\left\|\widehat{Z}_k-\widetilde{Z}H_k\right\|_F^2
  &=O_p\left(\frac{1}{m^2n^3}\right)\left(\left\|M_{ZE}\widehat{Z}_k\right\|_F^2
  +\left\|M_XM_{ZE}\widehat{Z}_k\right\|_F^2\right)\\
  &=O_p\left(\frac{1}{m^2n^3}\left\|M_{ZE}\widehat{Z}_k\right\|_F^2\right)
  =O_p\left(\frac{1}{\min\{m,n\}}\right),
  \end{aligned}
\end{equation}
where the last equality holds because the matrix $M_{ZE}$ does not include explanatory variables, and the Theorem 1 of \cite{bai2002determining} can be used. Based on equation (\ref{equ:error:z}), the Corollary 1 of \cite{bai2002determining} holds under Assumption \ref{assu1} and Assumption \ref{assu2}, which implies that
\begin{equation*}
    \lim_{\min\{m,n\}} P(\hat{q}=q)\to1.
\end{equation*}
Under the consistency of $\hat{q}$, we can assume a known $q$ without loss of generality. For the estimation error of $\hat{\gamma}_j$, we have
\begin{equation}\label{equ:error:gamma}
    \begin{aligned}
    \hat{\gamma}_j-H^{-1}\gamma_j=&H'\frac{1}{n}\sum_{i=1}^{n}z_i\varepsilon_{ij}
    -H'\frac{1}{n}Z'M_X\varepsilon_j+\frac{1}{n}\widehat{Z}'M_XZ\gamma_j\\
    &+\frac{1}{n}\widehat{Z}'\left(\widetilde{Z}H-\widehat{Z}\right)H^{-1}\gamma_j
    +\frac{1}{n}\left(\widehat{Z}-\widetilde{Z}H\right)'\varepsilon_j.
    \end{aligned}
\end{equation}
According to Bernstein's inequality of independent sub-exponential variables, for $j\leq m$, $k\leq q$ and $t\geq 0$, there are positive constants $c_1$ and $c_2$ satisfying
\begin{equation*}
    P\left(\left|\frac{1}{n}\sum_{i=1}^{n}z_{ik}\varepsilon_{ij}\right|>t\right)\leq 2\exp\left(-c_1n\min\left\{\frac{t^2}{c_2^2},\frac{t}{c_2}\right\}\right).
\end{equation*}
Thus, as $m\to\infty$, we have
\begin{equation*}
    \begin{aligned}
    &P\left(\max_{j\leq m,k\leq q}\left|\frac{1}{n}\sum_{i=1}^{n}z_{ik}\varepsilon_{ij}\right|>\frac{2c_2}{c_1}\sqrt{\frac{\log m}{n}}\right)\\
    \leq&\sum_{k=1}^{q}\sum_{j=1}^{m}P\left(\left|\frac{1}{n}
    \sum_{i=1}^{n}z_{ik}\varepsilon_{ij}\right|>\frac{2c_2}{c_1}\sqrt{\frac{\log m}{n}}\right)\leq 2qe^{-\log m}\to0.
    \end{aligned}
\end{equation*}
Equation (2.1) in \cite{fan2018large} indicates that  $\max\{\|H\|,\|H\|^{-1}\}=O_p(1)$. The first term of equation (\ref{equ:error:gamma}) is bounded as
\begin{equation*}
    \max_{j\leq m}\left\|\frac{1}{n}\sum_{i=1}^nz_i\varepsilon_{ij}\right\|=\max_{j\leq m}\sqrt{\sum_{k=1}^{q}\left(\frac{1}{n}\sum_{i=1}^{n}z_{ik}\varepsilon_{ij}\right)^2}
    =O_p\left(\sqrt{\frac{\log m}{n}}\right).
\end{equation*}
To control the second term of equation (\ref{equ:error:gamma}), consider
\begin{equation*}
    \varepsilon_j'M_X\varepsilon_j=_d
    \sigma_{\varepsilon,jj}\,\chi^2_j(p+1),
\end{equation*}
conditioning on $X$ with chi-square random variables $\chi^2_j(p+1)$ of arbitrary dependency. As $m\to\infty$, Markov's inequality implies that
\begin{equation*}
    \begin{aligned}
    &P\left(\max_{j\leq m}\frac{1}{n}\varepsilon_j'M_X\varepsilon_j>\frac{8\log m}{n}\right)\leq\sum_{j=1}^{m}P\left(\chi^2_j(p+1)>8\log m\right)\\
    =&\sum_{j=1}^{m}P\left(e^{\chi^2_j(p+1)/4}>e^{2\log m}\right)
    \leq\sum_{j=1}^{m}\mathbb{E}\left(e^{\chi^2_j(p+1)/4}\right)e^{-2\log m}\to0.
    \end{aligned}
\end{equation*}
The third term of equation (\ref{equ:error:gamma}) can be also bounded by chi-square random variables. For the last two terms of equation (\ref{equ:error:gamma}), denote $\Delta z_i=\hat{z}_i-H'\tilde{z_i}$, and for each $k\leq q$, we have
\begin{equation*}
  \begin{aligned}
  &\frac{1}{n}\sum_{i=1}^{n}\hat{z}_{ik}\,\Delta z_{ik}\leq
  \sqrt{\frac{1}{n}\sum_{i=1}^{n}
  \hat{z}_{ik}^2}\,\sqrt{\frac{1}{n}\sum_{i=1}^{n}\Delta z_{ik}^2}
  =O_p\left(\frac{1}{\min\{\sqrt{m},\sqrt{n}\}}\right),\\
  &\frac{1}{n}\sum_{i=1}^{n}\varepsilon_{ij}\,\Delta z_{ik}
  \leq\sqrt{\max_{j\leq m}\frac{1}{n}\sum_{i=1}^{n}
  \varepsilon_{ij}^2}\,\sqrt{\frac{1}{n}\sum_{i=1}^{n}\Delta z_{ik}^2}
  =O_p\left(\frac{1}{\min\{\sqrt{m},\sqrt{n}\}}\right).
  \end{aligned}
\end{equation*}
It follows that
\begin{equation*}
  \max_{j\leq m}\|\hat{\gamma}_j-H^{-1}\gamma_j\|
  =O_p\left(\sqrt{\frac{\log m}{n}}\right)+O_p\left(\frac{1}{\sqrt{m}}\right).
\end{equation*}
For the consistency of $\hat{\sigma}_{\varepsilon,jj}$, we have
\begin{equation*}
    \widehat{\Lambda}_{\varepsilon}=
    \frac{(\widetilde{Y}-\widehat{Z}\widehat{\Gamma})'(\widetilde{Y}-\widehat{Z}\widehat{\Gamma})}{n}
    =\frac{(\widetilde{E}+\widetilde{Z}\Gamma-\widehat{Z}\widehat{\Gamma})'
    (\widetilde{E}+\widetilde{Z}\Gamma-\widehat{Z}\widehat{\Gamma})}{n},
\end{equation*}
where each diagonal entry has the following decomposition
\begin{equation*}
    \hat{\sigma}_{\varepsilon,jj}=\frac{1}{n}\varepsilon_j'Q(\widetilde{X})\varepsilon_j+
    \frac{2}{n}\varepsilon_j'Q(\widetilde{X})(\widetilde{Z}\gamma_j-\widehat{Z}\hat{\gamma}_j)
    +\frac{1}{n}\|\widetilde{Z}\gamma_j-\widehat{Z}\hat{\gamma}_j\|^2,
\end{equation*}
and the estimation error is bounded as
\begin{equation}\label{equ:error:sigma}
    \begin{aligned}
    |\hat{\sigma}_{\varepsilon,jj}-\sigma_{\varepsilon,jj}|\leq
    &\left|\frac{1}{n}\varepsilon_j'Q(\widetilde{X})\varepsilon_j-\sigma_{\varepsilon,jj}\right|\\
    &+\frac{2}{n}\|\varepsilon_j\|\|\widetilde{Z}\gamma_j-\widehat{Z}\hat{\gamma}_j\|
    +\frac{1}{n}\|\widetilde{Z}\gamma_j-\widehat{Z}\hat{\gamma}_j\|^2.
    \end{aligned}
\end{equation}
Apply Bernstein's inequality again and the first term of equation (\ref{equ:error:sigma}) is controlled as
\begin{equation*}
    \begin{aligned}
    &\max_{j\leq m}\left|\frac{1}{n}\varepsilon_j'Q(\widetilde{X})\varepsilon_j-\sigma_{\varepsilon,jj}\right|\\
    \leq&\max_{j\leq m}\left|\frac{1}{n}\sum_{i=1}^{n}(\varepsilon_{ij}^2-\sigma_{\varepsilon,jj})\right|
    +\max_{j\leq m}\frac{1}{n}\varepsilon_j'\widetilde{X}(\widetilde{X}'\widetilde{X})^{-1}\widetilde{X}'\varepsilon_j
    =O_p\left(\sqrt{\frac{\log m}{n}}\right).
    \end{aligned}
\end{equation*}
For other terms in equation (\ref{equ:error:sigma}), we have
\begin{equation*}
    \max_{j\leq m}\frac{1}{\sqrt{n}}\|\varepsilon_j\|\leq\sqrt{\max_{j\leq m}
    \left|\frac{1}{n}\sum_{i=1}^{n}(\varepsilon_{ij}^2
    -\sigma_{\varepsilon,jj})\right|+\max_{j\leq m}\sigma_{\varepsilon,jj}}=O_p(1),
\end{equation*}
and
\begin{equation*}
    \begin{aligned}
    &\max_{j\leq m}\frac{1}{\sqrt{n}}\|\widetilde{Z}\gamma_j-\widehat{Z}\hat{\gamma}_j\|\\
    \leq&\max_{j\leq m}\frac{1}{\sqrt{n}}\left\|\widetilde{Z}H(\hat{\gamma}_j-H^{-1}\gamma_j)\right\|+
    \max_{j\leq m}\frac{1}{\sqrt{n}}\left\|(\widehat{Z}-\widetilde{Z}H)\hat{\gamma}_j\right\|\\
    \leq&\max_{j\leq m}\|\hat{\gamma}_j-H^{-1}\gamma_j\|\|H\|\sqrt{\frac{1}{n}\sum_{i=1}^{n}\|z_i\|^2}+
    \max_{j\leq m}\|\gamma_j\|\sqrt{\frac{1}{n}\sum_{i=1}^{n}\|\hat{z}_i-H'\tilde{z_i}\|^2}\\
    =&O_p\left(\sqrt{\frac{\log m}{n}}\right)+O_p\left(\frac{1}{\sqrt{m}}\right).
    \end{aligned}
\end{equation*}

\section{Proof of Theorem \ref{the:fdp}}
\label{app:the}
\setcounter{equation}{0}
\renewcommand{\theequation}{C.\arabic{equation}}

Both $\mathrm{FDP}_{\cdot}(t)$ and $\mathrm{POW}_{\cdot}(t)$ are bounded and uniformly integrable, so one only needs to prove the convergence in probability. For fixed $t\in(0,1]$, consider the following inequality
\begin{equation*}
    \begin{aligned}
    &\left|\mathrm{FDP}_{ora}(t)-\mathrm{FDP}_{adj}(t)\right|\\
    \leq&\left[\frac{R_{ora}(t)R_{adj}(t)}{m^2}\right]^{-1}
    \left[\frac{|V_{ora}(t)-V_{adj}(t)|\,V_{ora}(t)}{m_0^2}\right.\\
    &\left.+\frac{|S_{ora}(t)-S_{adj}(t)|\,V_{ora}(t)}{m_1m_0}+
    \frac{|V_{ora}(t)-V_{adj}(t)|\,R_{ora}(t)}{m_0m}\right],
    \end{aligned}
\end{equation*}
where $|V_{ora}(t)-V_{adj}(t)|/m_0$ and $|S_{ora}(t)-S_{adj}(t)|/m_1$ are two important terms, for example, $|S_{ora}(t)-S_{adj}(t)|/m_1$ directly implies the convergence of power. Let $A_{j}=\{\omega:p_{\,ora,j}\leq t,\ p_{\,adj,j}> t\}$ and $B_{j}=\{\omega:p_{\,ora,j}> t,\ p_{\,adj,j}\leq t\}$ be discrepant events. For arbitrary $\epsilon>0$ and $j\in\mathcal{I}_0$, we have
\begin{equation*}
    \begin{aligned}
    P(A_j)&=P(t+p_{\,ora,j}-p_{\,adj,j}<p_{\,ora,j}\leq t)\\
    &\leq P(t-\epsilon\leq p_{\,ora,j}\leq t)+P(|p_{\,ora,j}-p_{\,adj,j}|>\epsilon)\\
    &\leq \epsilon+P(\max_{j\in\mathcal{I}_0}|p_{\,ora,j}-p_{\,adj,j}|>\epsilon).
    \end{aligned}
\end{equation*}
Similarly, $P(B_j)$ has the same upper bound, so the difference of the number of false discoveries is bounded as
\begin{equation*}
    \begin{aligned}
    \mathbb{E}\left|\frac{V_{ora}(t)-V_{adj}(t)}{m_0}\right|
    &\leq\frac{1}{m_0}\sum_{j\in\mathcal{I}_0}\left[P(A_{j})+P(B_{j})\right]\\
    &\leq2\epsilon+2P(\max_{j\in\mathcal{I}_0}|p_{\,ora,j}-p_{\,adj,j}|>\epsilon).
    \end{aligned}
\end{equation*}
For $j\in\mathcal{I}_1$, we also have
\begin{equation*}
    \begin{aligned}
    &\mathbb{E}\left|\frac{S_{ora}(t)-S_{adj}(t)}{m_1}\right|
    \leq\frac{1}{m_1}\sum_{j\in\mathcal{I}_1}\left[P(A_{j})+P(B_{j})\right]\\
    \leq&\max_{j\in\mathcal{I}_1}P(t-\epsilon\leq p_{ora,j}\leq t+\epsilon) + 2P(\max_{j\in\mathcal{I}_1}|p_{\,ora,j}-p_{\,adj,j}|>\epsilon)\\
    \leq&2\epsilon+2P(\max_{j\in\mathcal{I}_1}|p_{\,ora,j}-p_{\,adj,j}|>\epsilon),
    \end{aligned}
\end{equation*}
where the last inequality holds because $\min_{j\in\mathcal{I}_1}|\mu_j|\asymp\sqrt{n}$. If the consistency of p-values is obtained, then as $\min\{m,n\}\to\infty$, we have
\begin{equation*}
    \begin{aligned}
    &\frac{1}{m_0}\left|V_{ora}(t)-V_{adj}(t)\right|\to_{L_1}0,\\
    &\frac{1}{m_1}\left|S_{ora}(t)-S_{adj}(t)\right|\to_{L_1}0.
    \end{aligned}
\end{equation*}
In addition, Appendix \ref{app:ora} shows that the rejection rate of the oracle procedure has a lower bound almost surely
\begin{equation*}
    \frac{R_{ora}(t)}{m}=\frac{m_0}{m}\frac{V_{ora}(t)}{m_0}+\frac{m_1}{m}\frac{S_{ora}(t)}{m_1}\geq t,
\end{equation*}
which means that
\begin{equation*}
    \begin{aligned}
    \left[\frac{R_{ora}(t)R_{adj}(t)}{m^2}\right]^{-1}
    &\leq\left(\frac{m}{R_{ora}(t)}\right)^2
    \left[1-\frac{m}{R_{ora}(t)}\left|\frac{R_{ora}(t)-R_{adj}(t)}{m}\right|\right]^{-1},
    \end{aligned}
\end{equation*}
which is bounded in probability. Therefore, Theorem \ref{the:fdp} is proved if we can show
\begin{equation*}
    \max_{j\in\mathcal{I}_0}|p_{\,ora,j}-p_{\,adj,j}|\to_p0\quad\text{and}\quad \max_{j\in\mathcal{I}_1}|p_{\,ora,j}-p_{\,adj,j}|\to_p0.
\end{equation*}
In this respect, we prove the part of the FAT-DW in advance and consider the following inequality
\begin{equation*}
    \left|p_{ora,j}-p_{adj,j}\right|=
    2\int_{\min\{-|t_{ora,j}|,\,-|t_{adj,j}|\}}^{\max\{-|t_{ora,j}|,\,-|t_{adj,j}|\}}\,\frac{1}{\sqrt{2\pi}}e^{-x^2/2}dx
    \leq\sqrt{\frac{2}{\pi}}\left|t_{ora,j}-t_{adj,j}\right|,
\end{equation*}
where the difference of t-scores can be expressed as
\begin{equation}\label{equ:error:t}
    \begin{aligned}
    t_{ora,j}-t_{adj,j}=&\frac{1}{\sqrt{\sigma_{\varepsilon,jj}\hat{\sigma}_{\varepsilon,jj}}}
    \left\{\frac{\hat{\sigma}_{\varepsilon,jj}-\sigma_{\varepsilon,jj}}
    {\sqrt{\hat{\sigma}_{\varepsilon,jj}}+\sqrt{\sigma_{\varepsilon,jj}}}\,
    (c_n\alpha_j+\eta_j)\right.\\
    &+\sqrt{\sigma_{\varepsilon,jj}}(H\hat{\gamma}_j-\gamma_j)'\zeta\\
    &+c_n\sqrt{\sigma_{\varepsilon,jj}}\,\hat{\gamma}_j'(\widehat{\Gamma} Q(\mathbbm{1}_m)\widehat{\Gamma}')^{-1}H^{-1}\Gamma Q(\mathbbm{1}_m)\alpha\\
    &+c_n\sqrt{\sigma_{\varepsilon,jj}}\,\hat{\gamma}_j'(\widehat{\Gamma} Q(\mathbbm{1}_m)\widehat{\Gamma}')^{-1}(\widehat{\Gamma}-H^{-1}\Gamma) Q(\mathbbm{1}_m)\alpha\\
    &+\sqrt{\sigma_{\varepsilon,jj}}\,\hat{\gamma}_j'(\widehat{\Gamma} Q(\mathbbm{1}_m)\widehat{\Gamma}')^{-1}H^{-1}\Gamma Q(\mathbbm{1}_m)\eta\\
    &+\sqrt{\sigma_{\varepsilon,jj}}\,\hat{\gamma}_j'(\widehat{\Gamma} Q(\mathbbm{1}_m)\widehat{\Gamma}')^{-1}(\widehat{\Gamma}-H^{-1}\Gamma) Q(\mathbbm{1}_m)\eta\\
    &\left.+\sqrt{\sigma_{\varepsilon,jj}}\,\hat{\gamma}_j'(\widehat{\Gamma} Q(\mathbbm{1}_m)\widehat{\Gamma}')^{-1}
    \widehat{\Gamma}Q(\mathbbm{1}_m)(\Gamma-H\widehat{\Gamma})'\zeta\,\right\}.
    \end{aligned}
\end{equation}
Since $(\widehat{\Gamma} Q(\mathbbm{1}_m)\widehat{\Gamma}')^{-1}$ repeatedly appears on the right-hand side of equation (\ref{equ:error:t}), we first discuss its upper bound. Consider the following decomposition
\begin{equation*}
    \begin{aligned}
    \frac{1}{m}H\widehat{\Gamma} Q(\mathbbm{1}_m)\widehat{\Gamma}'H'=&\left[\frac{1}{m}H\widehat{\Gamma} Q(\mathbbm{1}_m)\widehat{\Gamma}'H'-\frac{1}{m}\Gamma Q(\mathbbm{1}_m)\Gamma'\right]\\
    &+\left(\frac{1}{m}\Gamma Q(\mathbbm{1}_m)\Gamma'-\Sigma_{\gamma}\right)+\Sigma_{\gamma},
    \end{aligned}
\end{equation*}
which converges to $\Sigma_{\gamma}$ in probability by Appendix \ref{app:pca} and condition (C\ref{con:4}), so through continuous mapping theorem, we have
\begin{equation*}
    \left\|\left[\frac{1}{m}H\widehat{\Gamma} Q(\mathbbm{1}_m)\widehat{\Gamma}'H'\right]^{-1}\right\|
    =\lambda_{min}^{-1}\left[\frac{1}{m}H\widehat{\Gamma} Q(\mathbbm{1}_m)\widehat{\Gamma}'H'\right]
    =O_p\left(\lambda_{min}^{-1}\left(\Sigma_{\gamma}\right)\right).
\end{equation*}
Next, we obtain
\begin{equation*}
    \begin{aligned}
    &\max_{j\leq m}\,c_n\sqrt{\sigma_{\varepsilon,jj}}\,\left|\hat{\gamma}_j'(\widehat{\Gamma} Q(\mathbbm{1}_m)\widehat{\Gamma}')^{-1}(\widehat{\Gamma}-H^{-1}\Gamma) Q(\mathbbm{1}_m)\alpha\right|\\
    =&O_p\left(\left\|\frac{H'}{m\sqrt{n}}
    \sum_{i=1}^{n}z_i\sum_{j=1}^{m}\varepsilon_{ij}(\alpha_j-\bar{\alpha})\right\|\right)
    +O_p\left(\frac{\sqrt{n}}{\min\{\sqrt{m},\sqrt{n}\}}\frac{\|\Gamma Q(\mathbbm{1}_m)\alpha\|}{m}\right)\\
    =&O_p(m^{h-1})+O_p\left(\max\left\{\frac{\sqrt{n}}{m^{3/2}},\frac{1}{m}\right\}\|\Gamma Q(\mathbbm{1}_m)\alpha\|\right),
    \end{aligned}
\end{equation*}
where $\bar{\alpha}=\sum_j\alpha_j/m$ and the last equality holds because the variance of $\sum_{j=1}^{m}\varepsilon_{ij}(\alpha_j-\bar{\alpha})/m$ is bounded as
\begin{equation*}
    \frac{\alpha' Q(\mathbbm{1}_m)\Sigma_{\varepsilon}Q(\mathbbm{1}_m)\alpha}{m^2}\leq
    \frac{\|\Sigma_{\varepsilon}\|\|Q(\mathbbm{1}_m)\alpha\|^2}{m^2}=O\left(m^{2(h-1)}\right).
\end{equation*}
For other terms in equation (\ref{equ:error:t}), denote $\widetilde{\Gamma}=\Gamma Q(\mathbbm{1}_m)=(\tilde{\gamma}_1',\ldots,\tilde{\gamma}_q')'$. Let $\epsilon_i$ and $\tilde{\epsilon}_i$ be independent random variables from standard normal distribution. Then $\|\eta\|^2$ and $\|\widetilde{\Gamma}\eta\|^2$ can be represented as sum of independent random variables as
\begin{equation*}
    \begin{aligned}
    &\frac{1}{m}\|\eta\|^2=\frac{1}{m}\sum_{j=1}^{m}\sigma_{\varepsilon,jj}+
    \frac{1}{m}\sum_{j=1}^{m}\lambda_j(\epsilon_j^2-1)=O_p(1),\\
    &\frac{1}{m}\|\widetilde{\Gamma}\eta\|^2=\sum_{k=1}^{q}\tilde{\lambda}_k\,\tilde{\epsilon}_k^2=O_p(1),
    \end{aligned}
\end{equation*}
where $\lambda_j$ is the $j$th eigenvalue of $\Sigma_{\varepsilon}$ and $\tilde{\lambda}_k$ is the $k$th eigenvalue of $\widetilde{\Gamma}\Sigma_{\varepsilon}\widetilde{\Gamma}'/m$ with an upper bound as below
\begin{equation*}
    \begin{aligned}
    \max_{k\leq q}|\tilde{\lambda}_k|&\leq
    tr\left\{\frac{1}{m}\widetilde{\Gamma}\Sigma_{\varepsilon}\widetilde{\Gamma}\right\}
    =\sum_{k=1}^{q}\frac{1}{m}\tilde{\gamma}_k'\Sigma_{\varepsilon}\tilde{\gamma}_k=O\left(\lambda_{max}(\Sigma_{\varepsilon})\right).
    \end{aligned}
\end{equation*}
It follows that
\begin{equation*}
    \begin{aligned}
    &\max_{j\leq m}\sqrt{\sigma_{\varepsilon,jj}}\,\left|\hat{\gamma}_j'(\widehat{\Gamma} Q(\mathbbm{1}_m)\widehat{\Gamma}')^{-1}H^{-1}\Gamma Q(\mathbbm{1}_m)\eta\right|\\
    =&O_p\left(\frac{1}{m}\|\widetilde{\Gamma}\eta\|\right)=O_p\left(\frac{1}{\sqrt{m}}\right),
    \end{aligned}
\end{equation*}
and
\begin{equation*}
    \begin{aligned}
    &\max_{j\leq m}\sqrt{\sigma_{\varepsilon,jj}}\,\left|\hat{\gamma}_j'(\widehat{\Gamma} Q(\mathbbm{1}_m)\widehat{\Gamma}')^{-1}(\widehat{\Gamma}-H^{-1}\Gamma) Q(\mathbbm{1}_m)\eta\right|\\
    =&O_p\left(\left\|\frac{1}{m}(\widehat{\Gamma}-H^{-1}\Gamma)Q(\mathbbm{1}_m)\eta\right\|\right)\\
    =&O_p\left(\frac{\|\eta\|}{m}
    \left\{tr\left\{\left[(\widehat{\Gamma}-H^{-1}\Gamma)'(\widehat{\Gamma}-H^{-1}\Gamma)\right]^2\right\}\right\}^{1/4}\right)\\
    =&O_p\left(\frac{1}{\sqrt{m}}\sqrt{\sum_{j=1}^{m}\|\hat{\gamma}_j-H^{-1}\gamma_j\|^2}\right)
    =O_p\left(\sqrt{\frac{\log m}{n}}\right)+O_p\left(\frac{1}{\sqrt{m}}\right),
    \end{aligned}
\end{equation*}
and
\begin{equation*}
    \begin{aligned}
    &\max_{j\leq m}\sqrt{\sigma_{\varepsilon,jj}}\,\left|\hat{\gamma}_j'
    (\widehat{\Gamma}Q(\mathbbm{1}_m)\widehat{\Gamma}')^{-1}
    \widehat{\Gamma}Q(\mathbbm{1}_m)(\Gamma-H\widehat{\Gamma})'\zeta\right|\\
    =&O_p\left(\left\|\frac{1}{m}\widehat{\Gamma}Q(\mathbbm{1}_m)
    (\Gamma-H\widehat{\Gamma})'\zeta\right\|\right)\\
    =&O_p\left(\frac{\|\zeta\|}{m}
    \left\{tr\left\{\left[(\Gamma-H^{-1}\widehat{\Gamma})Q(\mathbbm{1}_m)
    \widehat{\Gamma}'\widehat{\Gamma}Q(\mathbbm{1}_m)
    (\Gamma-H^{-1}\widehat{\Gamma})'\right]^2\right\}\right\}^{1/4}\right)\\
    =&O_p\left(\max_{j\leq m}\|\hat{\gamma}_j-H^{-1}\gamma_j\|\right)=O_p\left(\sqrt{\frac{\log m}{n}}\right)+O_p\left(\frac{1}{\sqrt{m}}\right).
    \end{aligned}
\end{equation*}
In summary, since $\max_{j\leq m}|\eta_j|=O_p(\sqrt{\log m})$ and $\alpha_j\equiv0$ for $j\in\mathcal{I}_0$,
the consistency under true nulls is deduced as
\begin{equation*}
    \begin{aligned}
    \max_{j\in\mathcal{I}_0}|t_{ora,j}-t_{adj,j}|
    =&O_p(m^{h-1})+O_p\left(\frac{\log m}{\sqrt{n}}\right)\\
    &+O_p\left(\sqrt{\frac{\log m}{m}}\right)
    +O_p\left(\frac{\sqrt{n}}{m}\|\Gamma Q(\mathbbm{1}_m)\alpha\|\right).
    \end{aligned}
\end{equation*}
For the consistency under alternatives, recall $t_{ora,j}=\mu_j+\eta_j^*$ in Appendix \ref{app:ora}. For any $\delta\in(0,1)$, we have
\begin{equation*}
  P\left(|p_{ora,j}|>\delta\right)\leq P\left(|\eta_1^*|>\min_{j\in\mathcal{I}_1}|\mu_j|-\psi_{\delta/2}\right),
\end{equation*}
that exponentially decays to zero as $n\to\infty$. In addition, denote the difference of t-scores as $\Delta t_j=|t_{ora,j}-t_{adj,j}|$ that converges in the following form
\begin{equation*}
    \max_{j\in\mathcal{I}_1}\frac{\Delta t_j}{|\mu_j|}=o_p\left(1\right).
\end{equation*}
The adjusted p-value also decays to zero as
\begin{equation*}
    \begin{aligned}
    P\left(|p_{adj,j}|>\delta\right)&\leq P\left(|t_{ora,j}|<\Delta t_j+\psi_{\delta/2}\right)
    \leq P\left(|\eta_1^*|>|\mu_j|-\Delta t_j-\psi_{\delta/2}\right)\\
    &\leq P\left(|\eta_1^*|>\min_{j\in\mathcal{I}_1}|\mu_j|\left(1-\max_{j\in\mathcal{I}_1}\frac{\Delta t_j}{|\mu_j|}\right)-\psi_{\delta/2}\right)\\
    &\leq P\left(|\eta_1^*|>\min_{j\in\mathcal{I}_1}|\mu_j|\left(1-\delta\right)-\psi_{\delta/2}\right)
    + P\left(\max_{j\in\mathcal{I}_1}\frac{\Delta t_j}{|\mu_j|}>\delta\right).
    \end{aligned}
\end{equation*}
To prove the part of the AdaFAT, let $\hat{\mathcal{I}}_{ori,0}$ be the set that is not rejected by the original t-tests, and it suffices to focus on the error term in equation (\ref{equ:error:t}) as follows
\begin{equation*}
  c_n\sqrt{\sigma_{\varepsilon,jj}}\,\hat{\gamma}_j'(\widehat{\Gamma}_{ori} Q(\mathbbm{1}_{m_{ori,0}})\widehat{\Gamma}_{ori}')^{-1}H^{-1}\Gamma _{ori} Q(\mathbbm{1}_{m_{ori,0}})\alpha_{ori}=O_p\left(\sqrt{\frac{n\|\alpha_{ori}\|^2}{m_{ori,0}}}\right),
\end{equation*}
where $\Gamma_{ori}'$ and $\alpha_{ori}$ consist of the rows of $\Gamma'$ and $\alpha$ in $\hat{\mathcal{I}}_{ori,0}$ separately, and $m_{ori,0}$ is the cardinality of $\hat{\mathcal{I}}_{ori,0}$. Also, we have
\begin{equation*}
  \frac{\|\alpha_{ori}\|^2}{m_{ori,0}}=O_p\left(\frac{m_1-S_{ori}(t)}{m-R_{ori}(t)}\right),
\end{equation*}
and our problem reduces to show that
\begin{equation*}
  n\left[1-\frac{S_{ori}(t)}{m_1}\right]\to0\quad\text{and}\quad\limsup\frac{V_{ori}(t)}{m_0}<1,
\end{equation*}
in probability as $\min\{m,n\}\to\infty$. Without loss of generality, we assume $\tilde{\sigma}_{\varepsilon,jj}=\|\gamma_j\|^2+\sigma_{\varepsilon,jj}$ is known. According to Appendix \ref{app:ora}, the proof under alternatives holds for arbitrary correlations of p-values, se we have
\begin{equation*}
  n\mathbb{E}\left[1-\frac{S_{ori}(t)}{m_1}\right]=\frac{n}{m_1}\sum_{j\in\mathcal{I}_1}P(|t_{ori,j}|< \psi_{t/2})\leq nP\left(|\eta_1^*|>\min_{j\in\mathcal{I}_1}|\tilde{\mu}_j|-\psi_{t/2}\right),
\end{equation*}
where $\tilde{\mu}_j=c_n\alpha_j/\sqrt{\tilde{\sigma}_{\varepsilon,jj}}\asymp\sqrt{n}$. There exists a $\tilde{c}_t$ such that
\begin{equation*}
  nP\left(|\eta_1^*|>\min_{j\in\mathcal{I}_1}|\tilde{\mu}_j|-\psi_{t/2}\right)\leq 2ne^{-\tilde{c}_tn}\to0.
\end{equation*}
Next, conditioning on $\zeta$, we have $\max_{j\leq m}\gamma_j'\zeta/\sqrt{\sigma_{\varepsilon,jj}}=O_p(1)$ and
\begin{equation*}
  \frac{V_{ori}(t)}{m_0}\leq\frac{1}{m_0}\sum_{j\in\mathcal{I}_0}1_{\left\{\eta_j>\psi_{t/2}\sqrt{\tilde{\sigma}_{\varepsilon,jj}}-\gamma_j'\zeta\right\}}.
\end{equation*}
Using \cite{lyons1988strong}'s strong law of large numbers as in Appendix \ref{app:ora}, we obtain
\begin{equation*}
  \frac{1}{m_0}\sum_{j\in\mathcal{I}_0}1_{\left\{\eta_j^*>\psi_{t/2}-\max_{j\leq m}\gamma_j'\zeta/\sqrt{\sigma_{\varepsilon,jj}}\right\}}\to_{a.s.}1-\Phi\left(\psi_{t/2}-\max_{j\leq m}\frac{\gamma_j'\zeta}{\sqrt{\sigma_{\varepsilon,jj}}}\right)<1.
\end{equation*}
\end{appendices}

\bibliographystyle{AdaFAT}
\bibliography{AdaFAT}

\end{document}